

\documentstyle{amsppt}
\magnification=\magstep 1
\TagsOnRight
\NoBlackBoxes
\leftheadtext{V\.G\. Pestov}
\rightheadtext{Operator spaces}
\def\norm #1{{\left\Vert\,#1\,\right\Vert}}
\def\c #1{{C^\ast\langle\, #1\,\rangle}}

\def\C {{\Bbb C}}
\def\N{{\Bbb N}}
\def\e{{\epsilon}}

\def\MAX {{\text{\rm MAX\,}}}
\def\BALL {{\text{\rm BALL\,}}}

\def\Rep {{\text{\rm Rep}\,}}
\topmatter
\title
Operator spaces
and residually finite-dimensional $C^\ast$-algebras
\endtitle
\author
Vladimir G\. Pestov
\endauthor
\affil
Victoria University of Wellington, \\ P.O. Box 600, Wellington,  New
Zealand \\ \\
vladimir.pestov$\@$vuw.ac.nz
\endaffil
\abstract{For every operator space $X$ the
$C^\ast$-algebra containing it in a universal way is
residually finite-dimensional
(that is, has a separating family of
finite-dimensional representations). In particular,
the free $C^\ast$-algebra on any normed space so is.
This is an
extension of an earlier result by Goodearl and Menal, and our
short proof is based on a criterion
due to Exel and Loring.}

\endabstract
\endsubjclass
\endkeywords
\endtopmatter
\document

\smallpagebreak
\heading
1. Introduction
\endheading
\smallpagebreak

A $C^\ast$-algebra is called residually finite-dimensional
(RFD) if it has a separating family of finite-dimensional
representations. Clearly, every commutative $C^\ast$-algebra
is such; as shown by Choi \cite{5},
the $C^\ast$ group algebra on a free group with two
generators is RFD.
In 1990  Goodearl and Menal \cite{12} proved that
every $C^\ast$-algebra as an image of a
residually finite-dimensional one; they have shown that
the free $C^\ast$-algebra,
$\c\Gamma$,
on a set $\Gamma$ of free generators subject to the norm
restriction $\Vert x\Vert\leq 1,~x\in \Gamma$,
is residually finite-dimensional.
The most recent development is due to Exel and Loring \cite{8}
who proved that residual finite-dimensionality is
preserved by $C^\ast$-algebra coproducts.

The goal of this paper is to give an extension to the
aforementioned result by Goodearl and Menal
by invoking the concept of an operator space,
coming from quantized functional analysis
\cite{6}.

The construction of
the free $C^\ast$-algebra on a set is
the left adjoint functor
to the forgetful functor sending a $C^\ast$-algebra $A$
to the set of elements of its unit ball.
One can consider, however, ``less forgetful'' functors,
such as
$$A\mapsto \{\text{underlying
normed space of}~ A\},$$
or still a more subtle one
$$A\mapsto \{\text{underlying
operator space of}~ A\}$$
In both cases there exist
left adjoints, which we call the free $C^\ast$-algebra
on a normed (resp. operator) space. The former
construction is a
particular case of the latter one.
The Goodearl-Menal free $C^\ast$-algebra on a set $\Gamma$
is then exactly the free $C^\ast$-algebra on
the operator space of the form
$\MAX(l_1(\Gamma))$.

As the main result of this
work, we show that the free $C^\ast$-algebra
on any operator space is residually finite-dimensional.
Also, we demonstrate freeness of some subalgebras of
free $C^\ast$-algebras on operator spaces,
discuss free commutative $C^\ast$-algebras on normed spaces,
and conclude the paper with an open problem.

All $C^\ast$-algebras in this paper are assumed to be unital.
However, all results have their immedeate
non-unital analogs, obtained by proceeding to unitizations.

\smallpagebreak
\heading
2. Preliminaries
\endheading
\smallpagebreak

A {\it matrix norm} on a vector space $X$ is a
collection of norms $\norm{\cdot}_n, n\in\N$ on
the vector spaces $M_n(X)=X\otimes M_n(\C)$.
Every $C^\ast$-algebra carries a natural matrix norm.
Any linear map $f\colon X\to Y$ between two vector spaces
gives rise to a natural map $f_{(n)}\colon M_n(X)\to M_n(Y)$.
If $X$ and $Y$ carry matrix norms, then $f$ is said to be
{\it completely bounded} if
$\norm{f}_{cb}\overset def\to=
\sup\{\norm{f_{(n)}}_n\colon n\in\N\}
<+\infty$, a {\it complete contraction}
if $\norm{f}_{cb}\leq 1$,
and a {\it complete isometry} if every $f_{(n)}$ is an isometry.
A vector space $X$ together with a fixed matrix norm
is called an {\it (abstract) operator space}
if it is completely isometric to a subspace of
the algebra $\Cal B(H)$ of bounded operators on a
Hilbert space.
(Equivalently: to a subspace of a $C^\ast$-algebra.)
Operator spaces have been characterized by Ruan \cite{18}
as those spaces with matrix norms satisfying the three
conditions:
\smallskip
{(i)} $\norm{v\oplus\bold 0_m}_{n+m}=\norm v_n$,
\smallskip
{(ii)} $\max \{\norm{Bv}_n, \norm{vB}_n\}\leq
\norm B\norm v_n$,
\smallskip
{($L^\infty$)}
$\norm{v\oplus w}_{n+m}=
\max \{\norm v_n,\norm w_m\}$,
\smallskip
\noindent
for all $v\in M_n(X)$, $w\in M_m(X)$, and $B\in M_n$; here
$M_n=M_n(\C)$, $\bold 0_m$ is the null element of $M_m(X)$,
and

$$v\oplus w=
\left[\matrix
v & \bold 0 \\
\bold 0 & w
\endmatrix\right]$$

If $E$ is a normed space, then there is the maximal among all
matrix norms on $X$ making it into an operator space and such
that $\norm{}_1$ is the original norm on $X$.
The resulting operator space is denoted by
$\MAX(E)$. The correspondence
$E\mapsto\MAX(E)$ is functorial, and for
any bounded linear mapping
$f\colon E\to F$ one has $\norm f_{cb}=\norm f$ \cite{4, 7}.

\proclaim{Theorem 2.1  {\rm
(Noncommutative Hahn-Banach Theorem -- Arveson \cite{2},
Wittstock \cite{21}})}
Let $X$ be an operator space and $Y$ be an operator subspace.
Let $n\in\N$.
Then every completely bounded linear mapping
$f\colon Y\to M_n$ extends to a completely bounded mapping
$\tilde f\colon X\to M_n$ in such a way that
$\norm{f}_{cb}=\norm{\tilde f}_{cb}$.
\qed\endproclaim

The operator spaces play in noncommutative analysis
and geometry the same r\^ole as the normed spaces play
in commutative analysis and geometry \cite{6}.
For an account of
theory of operator spaces, including numerous examples,
the reader is referred to \cite{4, 6, 7, 18}.

If $A$ is a $C^\ast$-algebra,
then $\Rep(A,H)$ stands for the set of all
(degenerate and non-degenerate) representations of $A$ in $H$,
that is, $C^\ast$-algebra morphisms $A\to\Cal B(H)$,
endowed with the coarsest topology,
making all the mappings of the form
$$\Rep (A,H)\ni \pi\mapsto \pi(x)(\xi)\in H,~x\in A,~\xi\in H$$
continuous. Clearly, this topology is
inherited from
$C_s(A,\Cal B_s(H))$; here the subscript ``{\it s}''
stands for the topology of simple convergence, so that
$\Cal B_s(H)$ is the space $\Cal B(H)$ endowed with
the strong operator topology. The basic
neighbourhoods of an element $\pi\in\Rep(A,H)$ are of
the form
$$\Cal O_\pi[x_1,x_2,\dots,x_n;\xi_1,\xi_2,\dots,\xi_n;\e]$$
$$\overset def\to= \{\eta\in\Rep(A,H)\colon
\norm{\pi(x_i)(\xi_i)-\eta(x_i)(\xi_i)}<\e,~i=1,2,\dots, n\}$$
This topology was considered by Exel and Loring \cite{8},
and it is
finer than the Fell topology \cite{9}.
A representation $\pi\in\Rep(A,H)$ is termed
{\it finite-dimensional} if its essential space
(the closure of $\pi(A)(H)$) is finite-dimensional.
A $C^\ast$-algebra $A$ is called
{\it residually finite-dimensional}, or {\it RFD},
if it admits a family of finite-dimensional
representations which separate points.

\proclaim{Theorem 2.2 {\rm (Exel and Loring \cite{8})}}
A $C^\ast$-algebra $A$ is residually finite-dimensional
if and only if the set of finite-dimensional
representations is everywhere dense in $\Rep(A,H)$.
\qed\endproclaim

\smallpagebreak
\heading
3. Free $C^\ast$-algebras on operator spaces
\endheading
\smallpagebreak

The following result seems to be well known in the $C^\ast$-algebra
folklore.

\proclaim{Theorem 3.1}
Let $X$ be a normed space.
There exist a
$C^\ast$-algebra, $C^\ast\langle X\rangle $, and
an isometric embedding
$i_X\colon X\hookrightarrow C^\ast\langle X\rangle $ such that
$i_X(X)$ generates $C^\ast\langle X\rangle $ topologically,
and for every $C^\ast$-algebra $A$ and a
contractive mapping $f\colon X\to A$ there exists
a $C^\ast$-algebra morphism
$\hat f\colon C^\ast\langle X\rangle \to A$
with $f=\hat f\circ i_X$.
Such a pair $(C^\ast\langle X\rangle ,i)$ is essentially unique.
\qed\endproclaim

This result extends to the
class of operator spaces, and the argument remains pretty
standard. It follows the scheme of proving the
existence and uniqueness of
universal arrows in topological algebra and functional analysis,
invented by Kakutani for free topological groups \cite{14}
and used extensively since then in various situations
(cf. \cite{1, 10, 11, 13, 15-17}).

\proclaim{Theorem 3.2}
Let $X$ be an operator space.
There exist a $C^\ast$-algebra, $\c X$, and
a completely isometric embedding
$i_X\colon X\hookrightarrow \c X$ such that
$i_X(X)$ generates $\c X$ topologically,
and for every $C^\ast$-algebra $A$ and a completely
contractive mapping $f\colon X\to A$ there exists
a $C^\ast$-algebra morphism
$\hat f\colon \c X \to A$
with $f=\hat f\circ i_X$.
Such a pair $(\c X,i_X)$ is essentially unique.
\endproclaim

\demo{Proof}
Denote by $\frak F$ the class  of
all pairs $(A,j)$ where $A$ is a $C^\ast$-algebra
and $j: X\to L$ is a completely contractive linear
mapping such that
the image $j(X)$ topologically generates $A$.
By identifying the isomorphic pairs, one can assume that
$\frak F$ is a set.
Let $i_X$ stand for the diagonal product of mappings
$\Delta \{j\: (A,j)\in \frak F\}$, viewed as a mapping from $X$ to
the $C^\ast$-direct product
$\frak A = \prod_{(A,j)\in\frak F}A$.
Clearly, $i_X$ is correctly defined
and it is a complete contraction.
Denote by $\c X$ the $C^\ast$-subalgebra of the
$C^\ast$-algebra $\frak A$ generated by the image $i_X(X)$.
The universality and uniqueness of the pair
$(\c X, i_X)$ are checked immedeately.
To prove that $i_X$ is in fact a complete isometry, notice that
by the very definition of an operator space,
there is a pair $(A,j)\in\frak F$
such that $j\colon X\to A$ is a complete isometry.
\qed\enddemo

In category theoretic terms, the pair $(C^\ast\langle X\rangle ,i)$ is
the {\it universal arrow}
\cite{15}
to the forgetful functor, $\Cal F$,
from the category of $C^\ast$-algebras and $C^\ast$-algebra morphisms
to the category of operator spaces and complete contractions.
The above result shows that $\Cal F$ has left adjoint,
namely $\c\cdot$. This is another addition to the sphere of
practival applicability of the Blackadar's
general construction of
a $C^\ast$-algebra in terms of generators and relations \cite{3}.

Theorem 3.2 is indeed an extension of theorem 3.1.

\proclaim{Proposition 3.3}
Let $X$ be a normed space. Then the free $C^\ast$-algebra on
$X$ is canonically isomorphic to the free
$C^\ast$-algebra on the maximal operator space
$\MAX(X)$ associated to $X$.
\endproclaim

\demo{Proof}
Follows from the fact that the functor $X\mapsto \MAX(X)$
is left adjoint to the forgetful functor
from the category of operator spaces and complete contractions
to the category of normed spaces and contractions.
\qed\enddemo

Now we put the above construction in connection with the
concept of the free $C^\ast$-algebra on a set $\Gamma$ of free generators
\cite{3, 12}.

\proclaim{Proposition 3.4}
Let $\Gamma$ be a set. Then the free $C^\ast$-algebra on
$\Gamma$ is canonically isomorphic to the free
$C^\ast$-algebra on the normed space $l_1(\Gamma)$ or,
what is the same, the free $C^\ast$-algebra on the operator space
$\MAX(l_1(\Gamma))$.
\endproclaim

\demo{Proof}
The functor of the form $\Gamma\mapsto l_1(\Gamma)$ is
left adjoint to the functor sending a normed space to the set
of elements of its unit ball.
\qed\enddemo

\smallpagebreak
\heading
4. Representations of free $C^\ast$-algebras
\endheading
\smallpagebreak

The central result of our paper is the following.

\proclaim{Theorem 4.1}
The free $C^\ast$-algebra on an operator space,
$\c X$, is residually finite-dimensional.
\endproclaim

\demo{Proof}
Let $\pi\in \Rep(\c X,H)$ for some Hilbert space $H$;
in view of the Exel-Loring criterion 2.2, it suffices to find a
finite-dimensinal representation
in an arbitrary neighbourhood, $U$, of $\pi$.
Assume that
$U=\Cal O_\pi[f_1,f_2,\dots,f_n;\xi_1,\xi_2,\dots,\xi_n;\e]$
for some $n\in\N$, $f_1,f_2,\dots,f_n\in\c X$,
$\xi_1,\xi_2,\dots,\xi_n\in H$, and $\e>0$.

Fix for each $i$ a star-polynomial $p_i$ in
variables belonging to $X$
with $\norm{f_i-p_i}<\e/2$.
There exists a finite collection of elements
$\{x_1,x_2,\dots, x_m\}\subset X\cup X^\ast$ such that
each $p_i$ is a polynomial in variables
$x_1,x_2,\dots, x_m$. Denote by $k$ an upper bound of
the degrees of
$p_1,p_2,\dots, p_n$ relative to $\{x_1,x_2,\dots, x_m\}$.
Let $V$ be the subspace of $H$ spanned by all elements of
the form $\xi_1,\xi_2,\dots, \xi_n$ and
$\pi(x_{i_1})\pi(x_{i_2})\dots\pi(x_{i_l})(\xi_i)$,
where $i_1,i_2,\dots,i_l\in\{1,2,\dots,m\}$,
$l\leq k$, $i=1,2,\dots, n$.
Denote by $p_V$ the orthogonal projection from $H$ onto
the finite-dimensional subspace $V$.

For every $x\in X$ let $\eta(x)=p_V\pi(x)p_V$.
Since both the left and the right
multiplication by an idempotent element
in a $C^\ast$-algebra $A$ are complete contractions with respect to
the natural operator space structure on $A$,
the linear mapping $\eta\colon X\to \Cal B(H)$
is completely contractive.
By the universality of $\c X$, the mapping $\eta$ lifts to
a $C^\ast$-algebra morphism $\hat\eta\colon\c X\to\Cal B(H)$,
that is,  $\hat\eta\in\Rep(\c X,H)$.
It follows from the definition of $V$ that
for every $i=1,2,\dots,n$ one has
$\hat\eta(p_i)(\xi_i)=\pi(p_i)(\xi_i)$ and therefore
for each $i=1,2,\dots, n$
$$\norm{\pi(f_i)(\xi_i)- \hat\eta(f_i)(\xi_i)}\leq $$
$$\norm{\pi(f_i)(\xi_i)-\pi(p_i)(\xi_i)}+
\norm{\pi(p_i)(\xi_i) -\hat\eta(p_i)(\xi_i)}+
\norm{\hat\eta(p_i)(\xi_i)-\hat\eta(f_i)(\xi_i)} <$$
$$\e/2+0+\e/2=\e,$$
that is, $\hat\eta\in U$.

Finally, let $q$ be any star-polynomial with
variables from $X$. Since for any $x\in X$
one has $\hat\eta(x)(H)\subset V$ and
$(\hat\eta(x))^\ast(H)\subset V$,
one concludes that $\hat\eta(q)(H)\subset V$
as well. Since $\hat\eta$ is continuous
and polynomials of the form $q$ are dense in
$\c X$, the essential space of $\hat\eta$
is $V$, that is, the representation $\hat\eta$ is
finite-dimensional.
\qed\enddemo

\proclaim{Corollary 4.2}
The free $C^\ast$-algebra over a normed space $E$,
$C^\ast\langle E\rangle $, is residually finite-dimensional.
\qed\endproclaim

This was announced (without a proof) in our survey \cite{17}.

\proclaim{Corollary 4.3 {\rm (Goodearl and Menal \cite{12})}}
The free $C^\ast$-algebra over any set $\Gamma$,
$\c\Gamma$, is residually finite-dimensional.
\qed\endproclaim

\smallpagebreak
\heading
5. Free $C^\ast$-subalgebras
\endheading
\smallpagebreak

It is well known that if
$i\colon X\to Y$ is a monomorphism, then its extension
to universal objects, $\hat i$, need not
be a monomorphism anymore.
For example, the problem
of describing homeomorphic embeddings of topological spaces,
$i\colon X\to Y$, such that the monomorphism of
the free topological groups which extends $i$ is
topological, proved to be
fairly difficult
and was solved only recently \cite{19}.
On the contrary, the similar problem for free locally convex
spaces turned out to be readily amenable to the duality
techniques \cite{10, 11}.

The noncommutative Hahn-Banach theorem and our Theorem 4.1
enable us to solve the similar problem for free
$C^\ast$-algebras.

\proclaim{Theorem 5.1}
Let $X$ be an operator space and let $Y\hookrightarrow X$ be
an operator subspace.
Then the $C^\ast$-subalgebra of $\c X$,
generated by $Y$, is canonically isomorphic to the
free $C^\ast$-algebra on $Y$.
\endproclaim
\demo{Proof}
Denote by $i\colon Y\hookrightarrow X$ the complete isometry,
and by $\hat i\colon \c Y \to \c X$ its
extension to free $C^\ast$-algebras.
Being residually finite-dimensional,
the $C^\ast$-algebra $\c Y$ admits a monomorphism into
the $C^\ast$-direct product of a family of finite-dimensional
matrix algebras and therefore is isometric to
a $C^\ast$-subalgebra of such a $C^\ast$-direct product.
Now let $y\in \c Y$ and $\e>0$.
By the above said,
there exists an $n\in\N$ and a morphism
$f:\c Y\to M_n$ with
$$\norm{f(y)}_{M_n}>\norm y_{\c Y}-\e$$
By virtue of the noncommutative Hahn-Banach theorem 2.1,
the restriction $f_Y=f\vert_Y$, which is clearly
a complete contraction, extends to a complete contraction
$f_X\colon X\to M_n$.
The latter mapping lifts to a $C^\ast$-algebra morphism
$\tilde f\colon\c X\to M_n$, and
the following diagram commutes:

$$
\CD
 Y @> i_Y>> \c Y @> f>> M_n \\
@V i VV @V\bar i VV @| \\
X @>i_X>> \c X @> \tilde f>> M_n\\
\endCD
$$
Now it is clear that
$$\norm{\tilde f(\bar i(y))}_{M_n}=\norm{f(y)}_{M_n}>\norm y_{\c Y}-\e,$$
which means, in view of arbitrariness of an $\e>0$, that
$\norm{\bar i(y)}_{\c X}\geq \norm y_{\c Y}$,
and therefore
$\bar i:\c Y\to\c X$ is a $C^\ast$-algebra monomorphism, as desired.
\qed\enddemo

\smallpagebreak
\heading
6. The commutative case
\endheading
\smallpagebreak

The universal $C^\ast$-algebra on an operator space is essentially
non-commutative, and a sensible commutative analog of the construction
exists for normed spaces only: the construction
of the {\it free commutative $C^\ast$-algebra on
a normed space} $E$, which we denote by
$C^\ast_{com}\langle E\rangle$. It is the abelianization of the
algebra $\c E$.

\proclaim{Lemma 6.1}
The spectrum of the
free commutative $C^\ast$-algebra, $C^\ast_{com}\langle E\rangle$,
on an infinite-dimensional normed space $E$
is canonically homeomorphic to the closed unit ball
$\BALL(E_\sigma')$ of the weak dual space of $E$.
\endproclaim

\demo{Proof}
It is clear that the restriction map
$$\chi\mapsto\chi\vert_E\colon
\Sigma(C^\ast_{com}\langle E\rangle)\to \BALL(E_\sigma')$$
is one-to-one and continuous. Since the spectrum
$\Sigma(C^\ast_{com}\langle E\rangle)$ is compact,
the mapping is a homeomorphism.
\qed\enddemo

By $Q$ the Hilbert cube is denoted.
(Topologically, it is the countably
infinite Tychonoff power $I^{\aleph_0}$ of
the closed unit interval $I$.)

\proclaim{Theorem 6.2}
The free commutative $C^\ast$-algebra on any infinite-dimensional
separable normed space $E$ is isomorphic to
$C(Q)$.
\endproclaim

\demo{Proof}
According to Lemma 6.1,
$C^\ast_{com}\langle E\rangle\cong C(\BALL(E_\sigma'))$.
Since the unit ball of the weak dual space to a
separable normed space is homeomorphic to
the Hilbert cube $Q$
\cite{20, Lemma 3}, the statement follows.
\qed\enddemo

\smallpagebreak
\heading
7. Conclusion
\endheading
\smallpagebreak

The following
may be of some relevance to quantized functional analysis.

\definition{Problem 7.1}
Describe such pairs of operator spaces $X,Y$ that
their free $C^\ast$-algebras,
$\c X$ and $\c Y$, are isomorphic.
In particular, is it true or false that for any
infinite dimensional separable normed space $E$ the
algebra $\c E$
is isomorphic to the free $C^\ast$-algebra on
a countable set of free generators?
\enddefinition

A similar problem of classifying topological spaces
with isomorphic free topological groups
had lead to an independent, if somewhat scholastic,
problematics \cite{1, 13, 17}.

\smallpagebreak
\heading Acknowledgments
\endheading

This research was greatly stimulated by an
exchange of e-mail messages
with Terry Loring and N. Chris Phillips,
and I am thankful to both of them.
A Victoria University Small Research Grant
V212/451/RGNT/594/153 is acknowledged.

\smallpagebreak

\bye